\documentclass[reqno]{amsart}

\usepackage[english]{babel}
\usepackage{eucal,amsfonts,amstext}
\usepackage{amsmath}
\usepackage{amssymb}
\usepackage{amsthm}
\usepackage{setspace}
\usepackage{color}
\usepackage{mathtools}
\usepackage{enumerate}

\newtheorem{theorem}{Theorem}[section]
\newtheorem{corollary}{Corollary}[section]
\newtheorem{lem}{Lemma}[section]
\newtheorem{definition}{Definition}[section]
\newtheorem{exemple}{Example}[section]
\newtheorem{remark}{Remark}[section]
\newtheorem{proposition}{Proposition}[section]

\emergencystretch=\maxdimen
\begin{document}
	
\title[Stability and instability for discrete dynamics]{Weighted summability criteria for stability and instability of discrete dynamics with general growth rates}

\author[A. G\u ain\u a]{Ariana G\u ain\u a}
\address{A. G\u ain\u a, Department of Mathematics, Politehnica University of
Timi\c soara, Pia\c ta Victoriei 2, 300006 Timi\c soara, Romania}
\email{ariana.gaina@upt.ro}

\author[C. L\u azureanu]{Cristian L\u azureanu}
\address{C. L\u azureanu, Department of Mathematics, Politehnica University of
Timi\c soara, Pia\c ta Victoriei 2, 300006 Timi\c soara, Romania}
\address{Department of Mathematics, West University of Timi\c soara,
Bd. Vasile P\^arvan, No.4, 300223 Timi\c soara, Romania}
\email{cristian.lazureanu@upt.ro}

\author[N. Lupa]{Nicolae Lupa}
\address{N. Lupa, Department of Mathematics, Politehnica University of
Timi\c soara, Pia\c ta Victoriei 2, 300006 Timi\c soara, Romania}
\email{nicolae.lupa@upt.ro}

\begin{abstract}
Following the time rescaling introduced by Dragi\v{c}evi\'{c} and Silva in their paper \cite[Discrete Contin. Dyn. Syst. Ser. S 18 (2025), 3917--3944]{dragicevicsilva}, we establish several criteria for uniform $h$-stability and uniform $h$-instability of discrete evolution cocycles.
In particular,
these asymptotic properties are characterized by suitable weighted summability conditions under uniform $h$-boundedness and uniform $h$-decay assumptions, respectively. Explicit examples show
that the respective assumptions cannot
be omitted.
Our results cover a wide range of discrete growth rates $h$, including exponential and polynomial ones.
Furthermore, working in the general setting of discrete evolution cocycles, we provide a framework for discrete evolution families and discrete cocycles.
\end{abstract}
	


\maketitle
	
\section{Introduction}

In the asymptotic theory of both differential and difference equations, a central problem is to establish conditions ensuring exponential stability, instability, dichotomy, or trichotomy of solutions represented by evolution families \cite{Da,Hen,Po.2010}, cocycles over (semi)flows \cite{Ch.La.1999,Chow-Leiva-JDE-1995}, or, more generally, evolution cocycles over evolution (semi)flows \cite{hai4,PPP,SM-NA}.

A fundamental result in the stability theory of continuous-time evolution families  was obtained by Datko in  \cite[Theorem 1]{Da.1972}, where
he proved that a uniformly exponentially bounded evolution family $\left\{E(t,s)\right\}_{t\geq s\geq0}$
on a Banach space $W$ is uniformly exponentially stable (i.e., there exist constants $\nu>0$ and $N\geq 1$ such that $\| E(t,s)\|\leq Ne^{-\nu (t-s)}$ for all $t\geq s\geq 0$) if and only if, for some $p\in [1,\infty)$,
\[
\sup\limits_{s\geq0}\int_{s}^{\infty}\parallel E(\xi,s)w\parallel^{p}%
d\xi<\infty,\;w\in W.
\]

Since then, Datko's theorem has been widely recognized as one of the major achievements in control theory, leading to numerous generalizations and extensions in both continuous- and discrete-time settings; see, e.g., \cite{Be.Lu.Me.Si,davor,Lu.Me.Po,popaceausumegan,Pr.Me.1985,Rol,Sa.2010-1,Zabczyk95} for evolution families, \cite{dragicevic2018,DSSS2023,MSS-Stev2002,PPB,PPPCom,SMS2012,SS2004,SS2024,Sa.2010-1} for cocycles, and \cite{Gaina,hai4,pvhai,pvhaisi,MCT-ues,SM-NA} for evolution cocycles.

An advantage of the discrete-time case is that certain results obtained in this framework can be used to derive corresponding results in the continuous-time setting (see, for instance, \cite[Theorem 7.9, p. 213]{Ch.La.1999}, \cite{davor}, \cite{DSSS2023}, \cite[pp. 229-230]{Hen}, \cite{SS2004}, \cite{ZZ1}).

The first discrete-time version of the Datko theorem for nonautonomous linear difference equations was obtained by
Przy{\l}uski and Rolewicz in \cite{PR84}. More precisely,  they proved that a linear difference equation in a Banach space,
$$w_{n+1}=A_n w_n , \; n\geq n_0,$$ is uniformly exponentially stable (i.e., there exist constants  $\nu>0$ and $N\geq 1$ such that $\| E(m,n)\|\leq N e^{-\nu(m-n)}$ for all $m\geq n\geq n_0$ in $\mathbb{N}$) if and only if there exist $p\in[1,\infty)$ and $K>0$ such that
$$\sum\limits_{k=n}^\infty\|  E(k,n) \| ^p\leq K, \text{ for every } n\geq n_0,$$
where $E(m,n)=A_{m-1}\cdots A_n$ if $m>n$ in $\mathbb{N}$ and $E(n,n)=\mathrm{Id}$. For a comprehensive overview of this subject, we refer the reader to \cite{davor,SMS2012}.

The classical exponential behavior may appear too restrictive for the dynamics. It is therefore natural to consider more general types of asymptotic behavior, an approach that, to the best of our knowledge, was initiated by Martin \cite{Martin} and Muldowney  \cite{Mul}, who considered growth rates of the form  $e^{\int_0^t \rho(\xi) d\xi}$. We also refer to the work of Naulin and Pinto \cite{naulinpinto94,Pi.1994}, who studied the asymptotic behavior of evolution families with general growth rates.
More recently, Barreira and Valls developed a systematic framework for studying such phenomena in the context of nonuniform asymptotic behavior, which naturally arises when all Lyapunov exponents are either infinite or zero  \cite{Ba.Va.1} (see also the work of Bento and Silva \cite{BS2}). For more recent results in the theory of $h$-stability, $h$-instability, or $h$-dichotomy associated with general growth rates $h$,  we refer to \cite{borugamegantoth2022,dragicevic2026,dragicevicsilva,EPR,Ga.Ja2026,lupapopescu2023,silva2021}.

In our paper, we investigate uniform $h$-stability and uniform $h$-instability for discrete skew-product evolution semiflows, extending the classical framework based on the exponential growth rate $e^k$ to more general growth rates $h$.
In particular, we establish discrete Datko-type criteria characterizing these asymptotic properties for a broad class of discrete growth rates, including both exponential and polynomial ones. The question we address is how the summability condition in Datko-type criteria should depend on a prescribed growth rate $h$. We show that the sum appearing in our criteria involves the weight
$$\frac{h(k+1)-h(k)}{h(k)},\quad k\in\mathbb{N}.$$
For $h(k)=e^k$, the weight is constant, and hence the usual unweighted sum is recovered up to a constant factor. For $h(k)=k+1$, the weight is $\frac{1}{k+1}$. Both cases admit a common interpretation in terms of the time variable $\tau_k=\ln h(k)$:
the weight is comparable to  $\tau_{k+1}-\tau_k$ (see Remark \ref{remincrement}), and so each term in the sum is weighted according to the corresponding increment of the rescaled time.

Recently, several Datko-type characterizations of these types of behavior have been established for (discrete) evolution families in  \cite{RB,borugamegantoth2022,popa,glasnik}. However, the existing results  impose substantial restrictions on the growth rates, particularly by requiring them to belong to rather restrictive classes of function or sequence spaces.
Our approach, which is specific to discrete dynamics, is based on a time rescaling introduced by Dragi\v{c}evi\'{c} and Silva in \cite{dragicevicsilva}, and allows us to avoid the restrictions associated with these spaces.
The use of time rescaling has also proved instrumental in other contexts; see \cite{dragicevic2026,lupapopescu2023}.
	
This work is organized as
follows. The second section consists of three subsections. In the first subsection, we recall the notion of a discrete growth rate, introduce a continuous function associated with a given discrete growth rate, as defined by Dragi\v{c}evi\'{c} and Silva in \cite{dragicevicsilva}, provide illustrative examples, and establish several inequalities that play a crucial role in our analysis.
In the second subsection,
we consider the concept of a discrete skew-product evolution semiflow, which generalizes both discrete evolution families and discrete skew-product semiflows, thereby providing a unified framework for their study.
In the last subsection,
the continuous function introduced in subsection \ref{2.1} is used to construct a rescaled
discrete skew-product evolution semiflow (see Proposition \ref{prop_fih}).
The third section is devoted to characterizing uniform $h$-stability for discrete skew-product evolution semiflows and consists of three subsections. The first subsection focuses on uniform polynomial stability. More precisely, using a result due to Hai \cite{pvhai,pvhaisi}, we obtain a simple characterization of uniform polynomial stability (Proposition \ref{crit log ups}).
In the second subsection, following an approach developed by Dragi\v{c}evi\'{c} and Silva in \cite{dragicevicsilva}, we establish a connection between the uniform $h$-stability of the original discrete skew-product evolution semiflow and the uniform polynomial stability of the one introduced in Proposition \ref{prop_fih}, provided that $h$ satisfies a certain inequality.
In the third subsection, we first deduce a new criterion for uniform $h$-stability (Proposition \ref{critlog-uhs-discret}). This criterion plays a key role in proving the main result of this section, namely Theorem \ref{t1Zabczykhs}, which provides a  discrete Datko-type characterization of uniform $h$-stability.
To the best of our knowledge, this theorem is new even in the setting of discrete evolution families (Corollary \ref{evolf}),
and it naturally generalizes existing results in the literature.
The fourth section is devoted to the study of uniform $h$-instability. The constructions and techniques developed in the previous section are adapted to this setting. In the final section, we present two concrete examples of discrete skew-product evolution semiflows and apply Theorems \ref{t1Zabczykhs} and \ref{t1Zabczykhinst} to establish their uniform $h$-stability and uniform $h$-instability, respectively, for a suitable growth rate $h$.

\section{Preliminaries}
	
In this section, we recall the concepts of a discrete growth rate and a discrete skew-product evolution semiflow.
Furthermore, given a discrete growth rate and a discrete skew-product evolution semiflow, using the time rescaling introduced in \cite{dragicevicsilva}, we define a new discrete skew-product evolution semiflow that will be used in the subsequent analysis.

In what follows, $Y$ is a nonempty set, $W$ is a Banach space, and  $\mathcal{B}(W)$ denotes the Banach algebra of all bounded linear operators on $W$. To simplify
notation, the norm on $W$ and the operator norm on $\mathcal{B}(W)$ are both denoted by $\|\cdot\|$. Finally,  $\mathbb{N}=\{0,1,2,\ldots\}$ denotes the set of all non-negative integers, $\mathbb{N}^*=\mathbb{N}\setminus\{0\}$, and let
$$\Delta=\{(m,n)\in \mathbb{N}^2: m\geq n \}.$$

\subsection{Growth rates}\label{2.1}

We recall that a function $h:\mathbb{N}\to [1,\infty)$ is called a \emph{discrete growth rate} if it is strictly increasing  and satisfies
$$h(0)=1 \text{ and } \lim\limits_{n\to\infty}h(n)=\infty.$$

For instance, the functions $h(n)=e^n$ and $h(n)=n+1$, $n\in\mathbb{N}$, are both discrete growth rates, called the \emph{exponential growth rate} and the \emph{polynomial growth rate}, respectively.
While the theory of exponential behavior is well established and has been extensively developed over several decades (see, for instance, the monographs \cite{Ch.La.1999,Hen,Po.2010}), the corresponding theory for polynomial behavior is still evolving.
The notion of a polynomial growth rate was introduced independently by Barreira and Valls \cite{Ba.Va.4} and by Bento and Silva  \cite{Be.Si}. Another important growth rate considered in the literature is the \emph{logarithmic growth rate}, namely  $h(n)=\ln(n+e)$, $n\in\mathbb{N}$ (see \cite{silva2021}).

\medskip

Following the recent work of Dragi\v{c}evi\'{c} and Silva \cite{dragicevicsilva}, we associate with each discrete growth rate $h$  its piecewise linear interpolation
$\tilde{h}:[0,\infty)\to[1,\infty)$, defined by
$$\tilde{h}(t)=\begin{cases}
	h(n)& \text{ if } t=n\\
	h(n)+(t-n)(h(n+1)-h(n)) & \text{ if } n<t<n+1
\end{cases},\; n\in \mathbb{N}.$$

The function $\tilde{h}$ is continuous, strictly increasing and unbounded, with $\tilde{h}(0)=1$. Hence, it is a bijection from $[0,\infty)$ onto $[1,\infty)$, and its inverse
${\tilde{h}}^{-1}:[1,\infty)\to[0,\infty)$ is given by
$${\tilde{h}}^{-1}(t)=\begin{cases}
	n& \text{ if } t=	h(n)\\
	n+\frac{t-h(n)}{h(n+1)-h(n)} &\text{ if }	h(n)<t<	h(n+1)
\end{cases}, \;  n\in \mathbb{N}.$$
Evidently, ${\tilde{h}}^{-1}$ is continuous and strictly increasing.

In what follows, we establish some inequalities that play a crucial role in the sequel. Further inequalities for discrete growth rates can be found in \cite[Lemma 3.1]{silva2021}.
\begin{lem}\label{lema cu discret}
Let $h:\mathbb{N}\rightarrow [1,\infty)$ be a discrete growth rate such that there exists a constant $c\geq h(1) $ satisfying
\begin{equation}\label{romb}
\frac{h(n+1)}{h(n)}\leq c, \text{ for every } n\in \mathbb{N}.
\end{equation}

Then
\begin{equation}\label{steluta}
\frac{1}{c}\,n\leq h([\tilde{h}^{-1}(n)])\leq n, \; n\in \mathbb{N}^*,
\end{equation}
where $[a]$ denotes the integer part of a real number $a$, i.e., $[a]$ is the unique integer for which $a-1<[a]\leq a$.

Furthermore, for any fixed $\alpha>0$, the following inequality holds:
\begin{equation}\label{eq.halpha}
\sum\limits_{k=n}^{\infty} \frac{h(k+1)-h(k)}{h(k)^{\alpha+1}}\leq \frac{c^{\alpha+1}}{\alpha}\,\frac{1}{h(n)^{\alpha}},\; n\in\mathbb{N}.
\end{equation}
\end{lem}		

\begin{proof}
We first prove the inequalities in \eqref{steluta}. For any fixed $n\in\mathbb{N}^{*}$, one of the following alternatives holds.

{\it Case 1}. $n\in Range(h)$, that is, there exists $K=K(n)\in \mathbb{N}$ such that $h(K)=n$.
Then $$h([\tilde{h}^{-1}(n)])=h(K)=n.$$
Since $c\geq h(1)>1$, this yields that \eqref{steluta} holds.
	
{\it Case 2}. $n\notin Range(h)$. Then there exists  $K=K(n)\in \mathbb{N}$ such that

\begin{equation} \label{eq.ineg.nK}
h(K)<n<h(K+1).
\end{equation}
The definition of $\tilde{h}^{-1}$ implies $$\tilde{h}^{-1}(n)=K+\alpha_{K,\, n}, \text{where } \alpha_{K,\,n}=\frac{n-h(K)}{	h(K+1)-h(K)} \in (0,1).$$ Therefore $[\tilde{h}^{-1}(n)]=K,$ and hence $$h([\tilde{h}^{-1}(n)])=h(K)<n.$$
On the other hand, by  \eqref{romb} and \eqref{eq.ineg.nK}, we get $$h([\tilde{h}^{-1}(n)])=h(K)\geq \frac{1}{c}\, h(K+1)>\frac{1}{c}\,n.$$ This completes the proof of \eqref{steluta}.	

Let us now fix $\alpha>0$ and $m, n\in\mathbb{N}$ with $m>n$. Then
$$\int_{h(k)}^{h(k+1)} \frac{1}{t^{\alpha+1}}\,dt=\frac{1}{\alpha}\left( \frac{1}{h(k)^\alpha}-\frac{1}{h(k+1)^\alpha}\right),\; k= n, n+1,\ldots,m.$$
For $t\in [h(k),h(k+1)]$, as $h(k+1)\leq c\, h(k)$, we get
$$\frac{1}{t^{\alpha+1}}\geq \frac{1}{c^{\alpha+1}}\,\frac{1}{h(k)^{\alpha+1}}.$$
Thus,
$$\frac{1}{\alpha}\left( \frac{1}{h(k)^\alpha}-\frac{1}{h(k+1)^\alpha}\right)\geq \frac{1}{c^{\alpha+1}}\,\frac{h(k+1)-h(k)}{h(k)^{\alpha+1}},$$
or equivalently,
$$\frac{h(k+1)-h(k)}{h(k)^{\alpha+1}}\leq \frac{c^{\alpha+1}}{\alpha}\left( \frac{1}{h(k)^\alpha}-\frac{1}{h(k+1)^\alpha}\right).$$
Summing the above inequality over $k= n, \ldots, m$, yields
$$\sum\limits_{k=n}^{m} \frac{h(k+1)-h(k)}{h(k)^{\alpha+1}}\leq \frac{c^{\alpha+1}}{\alpha}\left(\frac{1}{h(n)^{\alpha}}-\frac{1}{h(m+1)^{\alpha}}\right),$$
and letting $m\to\infty$, we obtain \eqref{eq.halpha}.
\end{proof}

Observe that all three discrete growth rates mentioned above satisfy inequality \eqref{romb}. Furthermore, if $h$ is a discrete growth rate satisfying \eqref{romb} for some constant $c\geq h(1)$, then
$h(n)\leq c^n$ for every $n\in\mathbb{N}$.
Thus, $h$ grows no faster than an exponential function (see \cite{ZZ} for further results and discussions on sub-exponential growth rates). For order relations on discrete growth rates and their implications, we refer the reader to the recent work of Gallegos and Jara \cite{Ga.Ja2026}.
Nevertheless, the class of discrete growth rates satisfying \eqref{romb} contains a wide variety of examples.

\begin{exemple}\label{ex.gr}\rm
Let $s=(s_k)_{k\in\mathbb{N}}$ be a sequence of positive real numbers satisfying $$\sum\limits_{k=0}^\infty s_k=\infty,$$ and assume that there exists a constant $S>0$ such that $s_k\leq S$ for every $k\in \mathbb{N}$.
Then the function $h_s:\mathbb{N}\to [1,\infty)$, defined by
\begin{equation}\label{eq.hgr.s}
h_s(0)=1 \text{ and } h_s(n)=\prod_{k=0}^{n-1} \left(1+s_k\right) \text{ for } n\in\mathbb{N}^*,
\end{equation}
is a discrete growth rate satisfying
\begin{equation}\label{eq.S}
\frac{h_s(n+1)}{h_s(n)}=1+s_n\leq 1+S,\; n\in\mathbb{N},
\end{equation}
and hence \eqref{romb} holds for $c=1+S$.

In view of \eqref{eq.S}, it suffices to show that $\lim\limits_{n\to\infty}h_s(n)=\infty $. Indeed, using the inequality $$\ln (1+x)\geq \frac{x}{1+x},\quad x\geq 0,$$
together with the boundedness of $s$, we obtain
$$\ln h_s(n)=\sum\limits_{k=0}^{n-1} \ln \left(1+s_k\right)\geq \frac{1}{1+S}\sum\limits_{k=0}^{n-1} s_k\to \infty \text{ as } n\to \infty,$$
which proves the claim.
\end{exemple}

Significant results on the asymptotic behavior of discrete dynamical systems involving growth rates satisfying inequality \eqref{romb} can be found in \cite{BaDrZ,dragicevicsilva,silva2021}.

\subsection{Discrete evolution cocycles}

\begin{definition}\rm \label{def.semiflow}	
A mapping $\varphi:\Delta \times Y \to Y$ is called a \emph{discrete evolution semiflow} on $Y$ if the following conditions are satisfied:
	\begin{enumerate}[$(i)$]
		\item $\varphi(n,n,y)=y$, for all $(n,y)\in \mathbb{N} \times Y$;
		\item $\varphi(m,n,\varphi(n,n_0,y))=\varphi(m,n_0,y)$, for all  $m\geq n\geq n_0$ in $\mathbb{N}$ and $y\in Y$.
	\end{enumerate}
\end{definition}	

The concept of a discrete evolution semiflow extends that of a discrete semiflow. Indeed, if $\varphi_0:\mathbb{N}\times Y\to Y$ is a \emph{discrete semiflow} on $Y$, i.e.,
	\begin{enumerate}[$(i)$]
		\item $\varphi_0(0,y)=y$, for every $y\in Y$,
		\item $\varphi_0(m,\varphi_0(n,y))=\varphi_0(m+n,y)$, for all   $m, n\in\mathbb{N}$ and $y\in Y$,
	\end{enumerate}
then

\begin{equation}\label{sem}
\varphi(m,n,y):=\varphi_0(m-n,y), \; (m,n,y)\in\Delta\times Y,
\end{equation}
is a discrete evolution semiflow on $Y$.

\begin{definition}\rm \label{def.cocycle}	
An operator-valued function $\varPhi :\Delta \times Y \to \mathcal{B}(W)$  is called a {\it discrete evolution cocycle} over a discrete evolution semiflow $\varphi:\Delta \times Y \to Y$ if:
	\begin{enumerate}[$(i)$]
		\item $\varPhi(n,n,y)=\mathrm{Id}$ (the identity operator on $W$), for all $(n,y)\in \mathbb{N} \times Y$;
		\item $\varPhi(m,n,\varphi(n,n_0,y))\varPhi(n,n_0,y)=\varPhi(m,n_0,y)$, for all  $m\geq n\geq n_0$ in $\mathbb{N}$ and $y\in Y$.
	\end{enumerate}	

A discrete dynamics $\hat{\varphi}=\left(\hat{\varphi}_{m,\,n}\right)_{(m,\,n)\in\Delta}$, defined by
$$\hat{\varphi}_{m,\,n}:Y\times W \to Y\times W,\quad \hat{\varphi}_{m,\,n}(y,w)=(\varphi(m,n,y),\varPhi(m,n,y)w),$$ is called a  {\it discrete skew-product evolution semiflow} on $Y\times W$, where $Y$ is the base space and $W$ is the state space. For simplicity, we write $\hat{\varphi}=\left(\varphi,\varPhi\right)$.
\end{definition}

To the best of our knowledge, the concept of an evolution cocycle over an evolution semiflow was introduced by Megan, Stoica, and Buliga in \cite{MSB-Carpathian-2007}, in the continuous-time setting (see also \cite{hai4}, \cite{PPP}, and \cite{SM-NA}).

\emph{Discrete evolution cocycles extend the concepts of discrete evolution families and discrete cocycles}, thereby providing a unified framework for their study.

To this end, we first recall that a family $\left\{E(m,n)\right\}_{(m,n)\in\Delta}\subset \mathcal{B}(W)$ is called a \emph{discrete evolution family} on $W$ if
$$E(n,n)=\mathrm{Id} \text{ and } E(m,n)E(n,n_0)=E(m,n_0),\; m\geq n\geq n_0 \text{ in } \mathbb{N}.$$

It is straightforward to verify that if $\left\{E(m,n)\right\}_{(m,n)\in\Delta}$ is a discrete evolution family on $W$,  then the mapping
\begin{equation}\label{eq.ef}
\varPhi(m,n,y):=E(m,n),  \; (m,n,y)\in\Delta\times Y,
\end{equation}
is a discrete evolution cocycle over any discrete evolution semiflow on $Y$.

On the other hand,  a \emph{discrete cocycle} over a discrete semiflow $\varphi_0:\mathbb{N}\times Y\to Y$ is a mapping
$\varPhi_0:\mathbb{N}\times Y\to \mathcal{B}(W)$ satisfying
$$\varPhi_0(0,y)=\mathrm{Id} \text{ and } \varPhi_0(m,\varphi_0(n,y))\:\varPhi_0(n,y)=\varPhi_0(m+n,y),\; m,n\in\mathbb{N},\, y\in Y,$$
in which case the pair $\left(\varphi_0,\varPhi_0\right)$ is called a \emph{discrete skew-product semiflow} on $Y\times W$.

One can easily check that if $\left(\varphi_0,\varPhi_0\right)$ is  a {discrete skew-product semiflow} on $Y\times W$, then the mapping $\varPhi:\Delta\times Y\to \mathcal{B}(W)$, defined by $$\varPhi(m,n,y):=\varPhi_0(m-n,y),$$ is a discrete evolution cocycle over the discrete evolution semiflow $\varphi$ defined in \eqref{sem}, that is, $\left(\varphi,\varPhi\right)$ is  a {discrete skew-product evolution semiflow} on $Y\times W$.

In the classical theory of skew-product semiflows, the base space $Y$ is generally required to satisfy additional assumptions. For instance, in \cite{Ch.La.1999}, the base space is assumed to be a locally compact metric space. In contrast, our results do not impose any such restrictions on $Y$.

\subsection{Time rescaling}

The function $\tilde{h}^{-1}$ introduced in the first subsection allows us to define a new discrete skew-product evolution semiflow associated with a discrete growth rate $h$ and a discrete skew-product evolution semiflow. More precisely, the following result holds.

\begin{proposition}\label{prop_fih}
Let $h:\mathbb{N}\to [1,\infty)$ be a discrete growth rate, and let $\hat{\varphi}=(\varphi, \varPhi)$ be a discrete skew-product evolution semiflow on $Y\times W$. Define the mappings  $\varphi_{h}:\Delta\times Y\rightarrow Y$  and  $\varPhi_{h}:\Delta\times Y\rightarrow \mathcal{B}(W)$   by
$$\varphi_{h}(m,n,y)=\varphi([{\tilde{h}}^{-1}(m+1)],[{\tilde{h}}^{-1}(n+1)],y)$$ and
$$\varPhi_{h}(m,n,y)=\Phi([{\tilde{h}}^{-1}(m+1)], [{\tilde{h}}^{-1}(n+1)],y).$$
Then $\hat{\varphi}_h=\left(\varphi_h, \varPhi_h \right)$ is a discrete skew-product evolution semiflow on $Y\times W$.
\end{proposition}

\section{Uniform stability with growth rates}
	
The aim of this section is to establish a Datko-type criterion for uniform $h$-stability of discrete skew-product evolution semiflows. As particular cases, we obtain characterizations of uniform exponential stability and uniform polynomial stability.

Let $h:\mathbb{N}\to [1,\infty)$ be a discrete growth rate, and let  $\varPhi:\Delta\times Y\to\mathcal{B}(W)$ be a discrete evolution cocycle over a discrete evolution semiflow $\varphi:\Delta\times Y\to Y$.

\begin{definition}\label{defuhs}\rm
The discrete evolution cocycle $\varPhi$ (or the corresponding discrete skew-product evolution semiflow $\hat{\varphi}=(\varphi, \varPhi)$) is called:
\begin{itemize}
\item \emph{uniformly $h$-bounded} if there exist constants $\omega >0$ and $M\geq 1$ such that
\begin{equation}\label{eq.hg}
\|\varPhi(m,n,y)\|\leq M \left(\frac{h(m)}{h(n)}\right)^{\omega}, \text{ for all } (m,n,y)\in \Delta \times Y;
\end{equation}
\item \emph{uniformly $h$-stable} if there exist constants $\nu >0$ and $N\geq 1$ such that
\begin{equation}\label{eq.hs}
\|\varPhi(m,n,y)\|\leq N \left(\frac{h(m)}{h(n)}\right)^{-\nu}, \text{ for all }  (m,n,y)\in \Delta \times Y.
\end{equation}
\end{itemize}

In particular, if  \eqref{eq.hg} holds for $h(n)=e^n$ or $h(n)=n+1$, $n\in\mathbb{N}$, then  $\hat{\varphi}$ is called  \emph{uniformly exponentially bounded} or \emph{uniformly polynomially bounded}, respectively. Similarly, by considering these choices of $h$ in \eqref{eq.hs}, we obtain the classical concepts of \emph{uniform exponential stability} and  \emph{uniform polynomial stability}, respectively.
\end{definition}

\subsection{A criterion for uniform polynomial stability}

We begin by stating a discrete version of a result due to Hai \cite[Lemma 2.12]{pvhaisi} for continuous skew-product evolution semiflows (see also \cite[Lemma 3.1]{pvhai}).

\begin{lem}\label{lemHai}
Assume that  $\hat{\varphi}$ is uniformly  polynomially bounded. Then it is uniformly polynomially stable if and only if there exist constants
$k\in \mathbb{N}$, with $k\geq 2$, and $\theta\in (0,1)$  such that
	\begin{equation*}\label{eq.hai.p1}
		\|\varPhi(k \, n,n,y)\|\leq \theta, \text{ for all } (n,y)\in\mathbb{N}^* \times Y.
	\end{equation*}
\end{lem}

Using the above result, we establish a new characterization of uniform polynomial stability, which plays an essential role in our analysis.

\begin{proposition}\label{crit log ups}
Assume that  $\hat{\varphi}$ is uniformly  polynomially bounded. Then it is uniformly polynomially stable if  and only if there exists a constant  $K>0$ such that
\begin{equation}\label{eq.log.p}
||\varPhi(m,n,y)||\,\ln\frac{m+1}{\;n+1}\leq K, \text{ for all } (m,n,y)\in \Delta\times Y.
\end{equation}
\end{proposition}
\begin{proof}
{\it Necessity.} Assume that $\hat{\varphi}$ is uniformly polynomially stable with constants  $\nu>0$ and $N\geq 1$.
Then, for $(m,n,y)\in \Delta\times Y$,
\[
\|\varPhi(m,n,y)\|\ln\frac{m+1}{\;n+1} \leq  \frac{N}{\nu} \left(\frac{m+1}{\;n+1}\right)^{-\nu} \ln\left(\frac{m+1}{\;n+1}\right)^{\nu}\leq \frac{N}{\nu e},
\]
since $\frac{\ln{x}}{x}\leq \frac{1}{e}$ for every $x\geq 1$.
Hence,  \eqref{eq.log.p} holds for $K=\frac{N}{\nu e}$.

\medskip

{\it Sufficiency.} Assume that there exists a constant $K>0$ such that \eqref{eq.log.p} holds.
We claim that \emph{for every  $\theta\in (0,1)$  there exists  $k\in \mathbb{N}^*$ with $k\geq 2$  such that}
	\begin{equation*}\label{eq.hai.p}
		\|\varPhi(k \, n,n,y)\|\leq \theta, \text{ for all } (n,y)\in\mathbb{N}^* \times Y.
	\end{equation*}

\medskip

Indeed, let $\theta\in (0,1)$  and consider $k:=\left[2e^{\frac{K}{\theta}}\right]$. Since $k>2e^{\frac{K}{\theta}}-1$, or equivalently $\ln{\frac{k+1}{2}}>\frac{K}{\theta}$, for any
$(n,y)\in \mathbb{N}^*\times Y$, we have
$$\frac{K}{\theta} \,\|\varPhi(k \, n,n,y)\|\leq \ln{\frac{k+1}{2}}\, \|\varPhi(k \, n,n,y)\|\leq \ln{\frac{k\,n+1}{n+1} }\,\|\varPhi(k \, n,n,y)\|\leq K,$$
which proves the claim.

Finally,
Lemma \ref{lemHai} completes the proof.
\end{proof}

A result similar to Proposition \ref{crit log ups} for continuous-time evolution families was established in \cite[Theorem 3.1]{borugameganUVTln}.

\subsection{Uniform $h$-stability via uniform polynomial stability and time rescaling}
	
In the following result, based on the approach developed by Dragi\v{c}evi\'{c} and Silva in \cite{dragicevicsilva} for discrete evolution families, we establish a connection between the asymptotic behavior of the discrete skew-product evolution semiflows $\hat{\varphi}$ and $\hat{\varphi}_h$ considered in Proposition \ref{prop_fih}.

\begin{proposition}\label{legatura-uhs-ups-discret}
Let $h:\mathbb{N}\to [1,\infty)$ be a discrete growth rate satisfying the inequality
\eqref{romb} for some constant $c\geq h(1)$, and assume that $\hat{\varphi}=(\varphi, \varPhi)$ is a uniformly $h$-bounded discrete skew-product evolution semiflow on $Y\times W$. Then $\hat{\varphi}_{h}$ is uniformly polynomially bounded.

Furthermore,
$\hat{\varphi}$ is uniformly $h$-stable if and only if $\hat{\varphi}_{h}$ is uniformly polynomially stable.		
\end{proposition}
\begin{proof} We first prove that $\hat{\varphi}_{h}$ is uniformly polynomially bounded.
By Lemma \ref{lema cu discret}, for $(m,n,y)\in \Delta\times Y$,  we obtain
	\begin{align*}
		\|\varPhi_{h}(m,n,y)\|&=\|\varPhi([{\tilde{h}}^{-1}(m+1)], [{\tilde{h}}^{-1}(n+1)],y)\|\\
		&\leq M \left(\frac{h([{\tilde{h}}^{-1}(m+1)])}{h([{\tilde{h}}^{-1}(\;n+1)])}\right)^{\omega}\\
&\leq M c^{\omega} \left(\frac{m+1}{\;n+1}\right)^{\omega},
	\end{align*}
where $\omega >0$ and $M\geq 1$ satisfy \eqref{eq.hg}.
This yields that $\hat{\varphi}_{h}$ is uniformly polynomially bounded.	

Using the same arguments as before, one can easily prove that if $\hat{\varphi}$ is uniformly $h$-stable, then $\hat{\varphi}_{h}$ is uniformly polynomially stable.	
Conversely, assume that $\hat{\varphi}_{h}$ is uniformly polynomially stable with constants $\nu >0$ and $N\geq 1$.
Let $m,n\in\mathbb{N}$ with $m>n$. We distinguish the following two cases.

{\it Case 1}. There exist $i,j\in\mathbb{N}$ with $i>j$ such that
$$[\tilde{h}^{-1}(j+1)]\leq n< [\tilde{h}^{-1}(j+2)]\leq [\tilde{h}^{-1}(i+1)]\leq m < [\tilde{h}^{-1}(i+2)].$$

Using identity $(ii)$ in Definition \ref{def.cocycle}, the uniform  $h$-boundedness of $\hat{\varphi}$, inequality \eqref{steluta} and our assumption, we obtain
\begin{align*}
\| \varPhi(m,n,y)\|
&\leq\|\varPhi(m,[\tilde{h}^{-1}(i+1)],\varphi([\tilde{h}^{-1}(i+1)],n,y))\| \\
&\qquad  \| \varPhi([\tilde{h}^{-1}(i+1)],[\tilde{h}^{-1}(j+2)],\varphi([\tilde{h}^{-1}(j+2)],n,y))\|\\
&\qquad  \| \varPhi([\tilde{h}^{-1}(j+2)],n,y)\|\\
&\leq M^2 \left( \frac{h(m)}{h([\tilde{h}^{-1}(i+1)])} \right)^{\omega} \left( \frac{h([\tilde{h}^{-1}(j+2)])}{h(n)} \right)^{\omega}\\
&\qquad   \| \varPhi_h(i,j+1, \varphi([\tilde{h}^{-1}(j+2)],n,y))\|\\
&\leq M^2 c^\omega \left(\frac{i+2}{i+1}\right)^\omega c^\omega \left(\frac{j+2}{j+1}\right)^\omega N  \left(\frac{i+1}{j+2}\right)^{-\nu}\\
&\leq M^2 N (2c)^{2\omega}\left(\frac{j+2}{i+1}\right)^{\nu}\\
&\leq M^2 N (2c)^{2\omega}2^\nu \left(\frac{j+2}{i+2}\right)^{\nu},
\end{align*}
where we also used the inequalities $\frac{i+2}{i+1}\leq 2$ and  ${i+1}\geq  \frac{i+2}{2}$. Since
$$\frac{j+2}{i+2}=\frac{h(n)}{h(m)}\, \frac{h(m)}{i+2}\, \frac{j+1}{h(n)}\, \frac{j+2}{j+1}\leq 2c \left(\frac{h(m)}{h(n)}\right)^{-1}, $$
the above estimate yields
\begin{equation}\label{eq1.suf}
\| \varPhi(m,n,y)\|\leq M^2 N (2c)^{2\omega+\nu}2^\nu \left(\frac{h(m)}{h(n)}\right)^{-\nu}.
\end{equation}

{\it Case 2}. There exists $i\in\mathbb{N}$ such that
$$[\tilde{h}^{-1}(i+1)]\leq n<m\leq [\tilde{h}^{-1}(i+2)].$$

Then
\begin{align}\label{eq2.suf}
\| \varPhi(m,n,y)\|&\leq M \left( \frac{h(m)}{h(n)} \right)^{\omega}\notag\\
&= M \left( \frac{h(m)}{h(n)} \right)^{\omega+\nu}\left( \frac{h(m)}{h(n)} \right)^{-\nu}\notag\\
&\leq M c^{\omega+\nu} \left( \frac{i+2}{i+1} \right)^{\omega+\nu}\left( \frac{h(m)}{h(n)} \right)^{-\nu}\notag\\
&\leq M (2c)^{\omega+\nu}\left( \frac{h(m)}{h(n)} \right)^{-\nu}.
\end{align}
From \eqref{eq1.suf} and \eqref{eq2.suf}, there exists $N_1\geq 1$ such that
$$ \| \varPhi(m,n,y)\| \leq N_1 \left( \frac{h(m)}{h(n)} \right)^{-\nu},\; (m,n,y)\in\Delta\times Y,$$
and thus $\hat{\varphi}$ is uniformly $h$-stable. This completes the proof.
\end{proof}
	
\subsection{A discrete Datko-type criterion for uniform $h$-stability}

In this subsection, we assume that $\hat{\varphi}=(\varphi, \varPhi)$ is a discrete skew-product evolution semiflow on $Y\times W$  and  $h:\mathbb{N}\to [1,\infty)$ is a discrete growth rate satisfying the inequality \eqref{romb} for some constant $c\geq h(1)$.


In the following preliminary result, we extend Proposition \ref{crit log ups} to the setting of uniform $h$-stability.

\begin{proposition} \label{critlog-uhs-discret}
If $\hat{\varphi}$ is uniformly $h$-bounded, then it is uniformly $h$-stable if and only if there exists a constant $K>0$ such that
\begin{equation}\label{eq.ln.h}
\|\varPhi(m,n,y)\| \ln \frac{h(m)}{h(n)}\leq K, \text{ for all } (m,n,y)\in \Delta \times Y.
\end{equation}
\end{proposition}
	
\begin{proof}
The \emph{necessity} follows by arguments similar to those used in Proposition \ref{crit log ups}; therefore, we omit the proof.

To prove the \emph{sufficiency}, we use the equivalence between uniform $h$-stability of $\hat{\varphi}$ and uniform polynomial stability of $\hat{\varphi}_h$, established in Proposition \ref{legatura-uhs-ups-discret}, together with Proposition \ref{crit log ups}.
Indeed, assume that \eqref{eq.ln.h} holds for some constant  $K>0$, and  	
let $(m,n,y)\in\Delta\times Y$. By Lemma \ref{lema cu discret}, we have $$\frac{m+1}{\;n+1}\leq \frac{ c \: h ([\tilde h^{-1}(m+1)])}{\phantom{c\:} h ([\tilde h^{-1}(\;n+1)] )},$$ thus
\begin{equation}\label{eq.lnmn}
\ln\frac{m+1}{\;n+1}\leq  \ln c +\ln \frac{h ([\tilde h^{-1}(m+1)])}{h ([\tilde h^{-1}(\;n+1)] )}.
\end{equation}

If  $\ln \dfrac{h ([\tilde h^{-1}(m+1)])}{h ([\tilde h^{-1}(\;n+1)] )}	\geq 1$, then
$$\ln c\leq \ln c \cdot\ln \frac{h ([\tilde h^{-1}(m+1)])}{h ([\tilde h^{-1}(\;n+1)] )},$$
and so
\begin{equation}\label{eq.nou}
\ln\frac{m+1}{\;n+1}\leq (\ln c+1 )\ln \frac{h ([\tilde h^{-1}(m+1)])}{h ([\tilde h^{-1}(\;n+1)] )}.
\end{equation}
This yields
\begin{align*}
\|\varPhi_{h}&(m,n,y)\|\ln\frac{m+1}{\;n+1}\\
&= \|\varPhi([\tilde h^{-1}(m+1)], [\tilde h^{-1}(n+1)],y) \| \ln\frac{m+1}{\;n+1}\\
&\leq (\ln c+1 )\|\varPhi([\tilde h^{-1}(m+1)], [\tilde h^{-1}(n+1)],y) \| \ln \frac{h ([\tilde h^{-1}(m+1)])}{h ([\tilde h^{-1}(\;n+1)] )}\\
&\leq (\ln c+1 ) K.
\end{align*}
On the other hand, if $\ln \dfrac{h ([\tilde h^{-1}(m+1)])}{h ([\tilde h^{-1}(\;n+1)] )}<1$, or equivalently
$\dfrac{h ([\tilde h^{-1}(m+1)])}{h ([\tilde h^{-1}(\;n+1)] )}<e,$
then
\begin{align*}
\|\varPhi_{h}&(m,n,y)\|\ln\frac{m+1}{\;n+1}\\
&\leq  \|\varPhi([\tilde h^{-1}(m+1)], [\tilde h^{-1}(n+1)],y) \| \ln c  \\
&\quad +\|\varPhi([\tilde h^{-1}(m+1)], [\tilde h^{-1}(n+1)],y) \| \ln \frac{h ([\tilde h^{-1}(m+1)])}{h ([\tilde h^{-1}(\;n+1)] )}\\
&\leq M \left( \frac{h ([\tilde h^{-1}(m+1)])}{h ([\tilde h^{-1}(\;n+1)] )} \right)^\omega \ln c+K\\
&\leq M e^\omega \ln c +K.
\end{align*}
Therefore, the assumptions of Proposition \ref{crit log ups} are satisfied for $\hat{\varphi}_h$, and hence $\hat{\varphi}_h$ is uniformly polynomially stable. By Proposition \ref{legatura-uhs-ups-discret}, it follows that $\hat{\varphi}$ is uniformly $h$-stable.
\end{proof}	


The following theorem  provides a discrete Datko-type characterization of uniform $h$-stability for discrete skew-product evolution semiflows, considering a large class of discrete growth rates $h$.
To the best of our knowledge, this criterion is new even in the setting of discrete evolution families,
and it naturally generalizes existing results in the literature.

\begin{theorem}\label{t1Zabczykhs}
Let $h:\mathbb{N}\to [1,\infty)$ be a discrete growth rate such that there exists a constant $c\geq h(1)$ satisfying
\begin{equation}\label{eq.hth}
\frac{h(n+1)}{h(n)}\leq c, \text{ for every } n\in \mathbb{N},
\end{equation}
and assume that $\hat{\varphi}=(\varphi, \varPhi)$ is a uniformly $h$-bounded discrete skew-product evolution semiflow on $Y\times W$.
Then $\hat{\varphi}$ is uniformly $h$-stable if and only if
there exists a constant $D>0$ such that
\begin{equation}\label{eq.1}
\sum_{k=n}^{\infty}\frac{h(k+1)-h(k)}{h(k)}\,\|\varPhi(k,n,y)w\| \leq D \| w\|, \text{ for all } (n,y,w)\in \mathbb{N} \times Y\times W.
\end{equation}	
\end{theorem}
\begin{proof}
\emph{Necessity}.  Assume that $\hat{\varphi}$ is uniformly $h$-stable with constants $\nu>0$ and $N\geq 1$. Then, for $(n,y,w)\in\mathbb{N}\times Y\times W$, using inequalities  \eqref{eq.hs} and \eqref{eq.halpha}, we obtain
\begin{align*}
\sum_{k=n}^{\infty} \frac{h(k+1)-h(k)}{h(k)}\, \|\varPhi(k,n,y)w\|&\leq \sum_{k=n}^{\infty} \frac{h(k+1)-h(k)}{h(k)}\, N \left( \frac{h(k)}{h(n)}\right)^{-\nu}\|w\|\\
&=N h(n)^{\nu}\sum_{k=n}^{\infty}\frac{h(k+1)-h(k)}{h(k)^{\nu+1}}\|w\|\\
&\leq \frac{N \,c^{\nu+1}}{\nu}\|w\|,
\end{align*}
which proves the desired inequality.
		
\emph{Sufficiency}. For the reader's convenience, we divide the proof into two steps.

\medskip		
		
{\it Step 1.} \emph{We claim that there exists a constant $K\geq 1$ such that} $$\|\varPhi(m,n,y)\|\leq K,  \text{ for all } (m,n,y)\in \Delta \times Y.$$
			
\medskip

Indeed, let $(m,n,y,w)\in \Delta \times Y\times W$. If $h(m)\geq 2\,h(n)$, then $h(n)\leq \frac{h(m)}{2}<h(m)$ and thus there exists $p\in \{n, n+1,\ldots, m-1\}$ such that
$$h(p)\leq \frac{h(m)}{2}<h(p+1).$$
Since $h(p+1)\leq c\, h(p)$, the above inequality yields
$\frac{h(m)}{h(p)}<2\,c.$ Therefore,
$$ \frac{h(m)}{h(k)}\leq \frac{h(m)}{h(p)}<2\,c, \text{ for each fixed } k\in\{p, p+1, \ldots, m-1\}.$$
We now successively obtain
\begin{align*}
\|\varPhi&(m,n,y)w\|\\
&=\frac{1}{\sum\limits_{k=p}^{m-1} \frac{h(k+1)-h(k)}{h(k)}}\; \sum_{k=p}^{m-1} \frac{h(k+1)-h(k)}{h(k)}\,\|\varPhi(m,k,\varphi(k,n,y))\varPhi(k,n,y)w\|\\
&\leq  \frac{1}{\sum\limits_{k=p}^{m-1} \frac{h(k+1)-h(k)}{h(m)}}\; \sum\limits_{k=p}^{m-1} \frac{h(k+1)-h(k)}{h(k)}\, M \left(\frac{h(m)}{h(k)}\right)^{\omega}\|\varPhi(k,n,y)w\| \\
&\leq  \frac{M}{\frac{h(m)-h(p)}{h(m)}} \; \sum\limits_{k=p}^{m-1} \frac{h(k+1)-h(k)}{h(k)}\, (2\,c)^{\omega}\, \|\varPhi(k,n,y)w\| \\
&\leq  \frac{M (2\,c)^{\omega}}{ 1-\frac{h(p)}{h(m)} } \; \sum\limits_{k=n}^{\infty} \frac{h(k+1)-h(k)}{h(k)}\,  \|\varPhi(k,n,y)w\|,
\end{align*}
where $\omega>0$ and $M\geq 1$ satisfy \eqref{eq.hg}.
As $\frac{h(p)}{h(m)}\leq \frac{1}{2}$, this yields
\begin{equation}\label{eq.c1}
\|\varPhi(m,n,y)w\|\leq    M  D c^{\omega }2^{\omega+1 } \: \|w\|.
\end{equation}	
On the other hand, if  $h(m)<2\,h(n)$, we have
\begin{equation}\label{eq.c2}
\|\varPhi(m,n,y)w\|\leq M\left(\frac{h(m)}{h(n)}\right)^{\omega}\| w\| \leq 2^{\omega} M\|w\|.
\end{equation}
The claim follows from inequalities \eqref{eq.c1} and \eqref{eq.c2}.
			
\medskip
			
{\it Step 2.} \emph{We prove that}
\begin{equation}\label{eq.log}
\|\varPhi(m,n,y) \| \ln \frac{h(m)}{h(n)}\leq K D, \text{ for all } (m,n,y)\in\Delta\times Y.
\end{equation}

Indeed, let $(m,n,y,w)\in\Delta\times Y\times W$ with $m\geq n+1$. By Step 1 and \eqref{eq.1}, we obtain		
\begin{align*}
\|\varPhi(m,n,y)w\|\ln \frac{h(m)}{h(n)}&= \sum_{k=n}^{m-1} \|\varPhi(m,n,y)w\|\ln \frac{h(k+1)}{h(k)}\\
&= \sum_{k=n}^{m-1} \|\varPhi(m,k,\varphi(k,n,y))\varPhi(k,n,y)w\|\ln \frac{h(k+1)}{h(k)}\\
&\leq K \sum_{k=n}^{m-1}\|\varPhi(k,n,y)w\|\left( \frac{h(k+1)}{h(k)}-1\right)\\
&\leq K \sum_{k=n}^{\infty}  \frac{h(k+1)-h(k)}{h(k)}\,\|\varPhi(k,n,y)w\| \\
&\leq K D \|w\|,
\end{align*}
which proves \eqref{eq.log}. In the above estimates, we use the inequality $\ln x\leq x-1$ for $x=\frac{h(k+1)}{h(k)}\geq 1$.

\medskip
		
Finally, Proposition \ref{critlog-uhs-discret} completes the proof.			
\end{proof}
\begin{remark}\label{remincrement}
Assuming that \eqref{eq.hth} holds for some $c\geq h(1)$, we have
$$\tau_{k+1}-\tau_k\leq \frac{h(k+1)-h(k)}{h(k)} \leq c \left(\tau_{k+1}-\tau_k\right), \text{ where } \tau_k=\ln{h(k)},\; k\in\mathbb{N}.$$
This shows that the weight $w_k:=\frac{h(k+1)-h(k)}{h(k)}$ appearing in the estimate \eqref{eq.1} is closely related to the increments of the time rescaling $\tau_k=\ln h(k)$. This is natural in view of the identity
$$\left(\frac{h(m)}{h(n)}\right)^{-\nu}=e^{-\nu(\tau_m-\tau_n)},\quad (m,n)\in\Delta.$$
\end{remark}

As a consequence of Theorem \ref{t1Zabczykhs}, we deduce from \eqref{eq.ef} a natural characterization of uniform $h$-stability for discrete evolution families. We recall that a discrete evolution family  $\mathcal{E}=\left\{E(m,n)\right\}_{(m,n)\in\Delta}$ is called \emph{uniformly $h$-bounded}
if there exist constants $\omega>0$ and $M\geq 1$ such that
$$\| E(m,n) \|\leq M \left(\frac{h(m)}{h(n)}\right)^\omega,\; (m,n)\in\Delta.$$
If, instead, the inequality
$$\| E(m,n) \|\leq N \left(\frac{h(m)}{h(n)}\right)^{-\nu},\; (m,n)\in\Delta,$$
holds for some constants $\nu>0$ and $N\geq1$, then $\mathcal{E}$ is said to be \emph{uniformly $h$-stable}.

\begin{corollary}\label{evolf}
Let $h:\mathbb{N}\to [1,\infty)$ be a discrete growth rate such that \eqref{romb} holds for some $c\geq h(1)$, and assume that $\left\{E(m,n)\right\}_{(m,n)\in\Delta}$ is a uniformly $h$-bounded evolution family on $W$.
Then it is uniformly $h$-stable if and only if
there exists a constant $D>0$ such that
\begin{equation}\label{eq.1fe}
\sum_{k=n}^{\infty}\frac{h(k+1)-h(k)}{h(k)}\,\|E(k,n)w\| \leq D \| w\|, \text{ for all } (n,w)\in \mathbb{N}\times W.
\end{equation}
\end{corollary}

The following example shows that the uniform $h$-boundedness assumption in both Theorem \ref{t1Zabczykhs} and Corollary \ref{evolf} cannot be removed.
\begin{exemple}\label{ex.g}\rm
Let $$s_k=
\begin{cases}
1,& k=2j\\
\frac{1}{2^j},& k=2j+1
\end{cases},
\quad j\in\mathbb{N}.
$$
Since $0<s_k\leq 1$ for every $k\in\mathbb{N}$ and $\sum\limits_{k=0}^\infty s_k=\infty$,
from Example \ref{ex.gr}, we obtain that
the function $h_s:\mathbb{N}\to [1,\infty)$, defined by
$$ h_s(0)=1 \text{ and } h_s(n)=\prod_{k=0}^{n-1} \left(1+ s_k\right) \text{ for } n\in\mathbb{N}^*,$$
is a discrete growth rate satisfying inequality \eqref{romb} for $c=2$. Moreover,
$$\frac{h_s(k+1)-h_s(k)}{h_s(k)}=s_k,\quad k\in\mathbb{N}.$$

For an arbitrary base space $Y$ and a nontrivial state space $W$, we also consider the discrete evolution cocycle generated by the discrete evolution family
$$E(m,n)=\begin{cases} A_{m-1} \cdots A_n,& m>n\\
\mathrm{Id},& m=n\end{cases},
\text{ where } A_n=\begin{cases} 2^j \, \mathrm{Id},& n=2j\\ 0,& n=2j+1\end{cases}, \quad j\in\mathbb{N},$$
that is,
$$\varPhi(m,n,y):=E(m,n)=\begin{cases} 0,& m>n+1\\ A_n,& m=n+1 \\ \mathrm{Id},& m=n\end{cases},\quad (m,n,y)\in\Delta\times Y.$$

We claim that $\varPhi$ satisfies \eqref{eq.1} (and hence $\mathcal{E}=\{E(m,n)\}_{(m,n)\in\Delta}$ satisfies \eqref{eq.1fe}) but is not uniformly $h_s$-stable.
\end{exemple}

\begin{proof}[\unskip\nopunct]
Indeed, for $(n,y,w)\in \mathbb{N} \times Y\times W$,  we obtain
\begin{align*}
\sum_{k=n}^{\infty}\frac{h_s(k+1)-h_s(k)}{h_s(k)}\,\|\varPhi(k,n,y)w\|
&=\sum_{k=n}^{\infty} s_k \,\|E(k,n)w\|\\
&=s_n \,\|w\|+s_{n+1}\| A_n w\|\\
&\leq 2 \|w\|,
\end{align*}
thus \eqref{eq.1} holds for $D=2$. Since
$$\| \varPhi(2n+1,2n,y) \|=\| A_{2n} \|=2^n,\; (n,y)\in \mathbb{N}\times Y,$$
we obtain that $\varPhi$ is not uniformly $h_s$-stable (the same conclusion also holds for $\mathcal{E}$).
\end{proof}

By setting $h(n)=e^n$ in Theorem \ref{t1Zabczykhs}, we recover a discrete Datko-type characterization of uniform exponential stability (see \cite{hai4,MCT-ues,popaceausumegan,SM-NA} and the references therein).

\begin{corollary}\label{cexp}
A uniformly exponentially bounded discrete skew-product evolution semiflow $\hat{\varphi}=(\varphi, \varPhi)$ is uniformly exponentially stable if and only if there exists a constant $D>0$ such that
\begin{equation*}\label{eq.expD}
\sum_{k=n}^{\infty} \|\varPhi(k,n,y)w\| \leq D \|w\|, \text{ for all } (n,y,w)\in \mathbb{N} \times Y\times W.
\end{equation*}
\end{corollary}

In particular, from the above corollary, we obtain a Datko-type characterization of uniform exponential stability for discrete cocycles over discrete semiflows (for similar results, we refer to  \cite{dragicevic2018,DSSS2023,SS2004}). Let us note that, for cocycles, it is natural to focus on exponential growth rates, rather than on other types of asymptotic behavior.

We first recall that a \emph{discrete cocycle} $\varPhi_0:\mathbb{N}\times Y\to \mathcal{B}(W)$  (over a discrete semiflow $\varphi_0:\mathbb{N}\times Y\to Y$) is
called \emph{uniformly exponentially bounded} if there exist constants $\omega>0$ and $M\geq 1$ satisfying
$$\|\varPhi_0(n,y)\|\leq M e^{\omega n}, \text{ for all $(n,y)\in\mathbb{N}\times Y$.}$$
Replacing $\omega$ with $-\omega$ in the above inequality, we obtain the notion of \emph{uniform exponential stability} for discrete cocycles.

\begin{corollary}\label{cexpcs}
Assume that $\varPhi_0:\mathbb{N}\times Y\to\mathcal{B}(W)$ is a uniformly exponentially bounded discrete cocycle over a discrete semiflow $\varphi_0:\mathbb{N}\times Y\to Y$. Then $\varPhi_0$ is uniformly exponentially stable if and only if there exists a constant $D>0$ such that
\begin{equation*}\label{eq.expDcs}
\sum_{k=0}^{\infty} \|\varPhi_0(k,y)w\| \leq D \|w\|, \text{ for all } (y,w)\in Y\times W.
\end{equation*}
\end{corollary}

Considering $h(n)=n+1$ in Theorem \ref{t1Zabczykhs}, we obtain a discrete Datko-type criterion for uniform polynomial stability. Similar results have been obtained in \cite[Corollary 3.4]{pvhai}, \cite[Corollary 3.2]{pvhaisi}, and \cite[Theorem 2.1]{tomescu}, although different methods were used.
\begin{corollary}\label{cpol}
A uniformly polynomially bounded discrete skew-product evolution semiflow $\hat{\varphi}=(\varphi, \varPhi)$ is uniformly polynomially stable if and only if there exists a constant $D>0$ such that
\begin{equation*}\label{eq.polD}
\sum_{k=n}^{\infty} \frac{\|\varPhi(k,n,y)w\|}{k+1} \leq D \|w\|, \text{ for all } (n,y,w)\in \mathbb{N} \times Y\times W.
\end{equation*}
\end{corollary}

The above corollaries show that Theorem \ref{t1Zabczykhs} incorporates several important results from the literature. However, our approach differs from those used to obtain these results.

\section{Uniform instability with growth rates}

In this section, we develop a theory of uniform instability with general growth rates following the same approach as that used for the stability theory in the previous section.

We continue to assume that  $h:\mathbb{N}\to [1,\infty)$ is a discrete growth rate  and that  $\hat{\varphi}=(\varphi, \varPhi)$ is a discrete skew-product evolution semiflow on $Y\times W$. We emphasize that no invertibility assumption is imposed on the operators $\Phi(m,n,y)$.
	
First, we introduce the notions of \emph{uniform $h$-decay} and \emph{uniform $h$-instability}.
\begin{definition}\label{defuhis}\rm
We say that	$\hat{\varphi}$  has
		\emph{uniform $h$-decay} if there exist constants $\omega >0$ and $M\geq 1$ such that
\begin{equation}\label{hbd}
M \|\varPhi(m,n,y)w\| \geq   \left(\frac{h(m)}{h(n)}\right)^{-\omega}\| w \|, \text{ for all } (m,n,y,w)\in \Delta \times Y\times W.
\end{equation}
In particular, $\hat{\varphi}$ is called  \emph{uniformly $h$-unstable} if there exist constants $\nu >0$ and $N\geq 1$  such that
\begin{equation}\label{hins}
N \|\varPhi(m,n,y)w\| \geq   \left(\frac{h(m)}{h(n)}\right)^{\nu}\| w \|, \text{ for all } (m,n,y,w)\in \Delta \times Y\times W.
\end{equation}

By setting $h(n)=e^n$ and $h(n)=n+1$ for every $n\in\mathbb{N}$, we recover the notions of \emph{uniform exponential decay} (\emph{uniform exponential instability}) and \emph{uniform polynomial decay} (\emph{uniform polynomial instability}), respectively.
\end{definition}

In the following result, we establish a Hai-type criterion for uniform polynomial instability of discrete skew-product evolution semiflows.

\begin{lem}\label{crithaiupis}
Assume that $\hat{\varphi}$ has uniform polynomial decay. Then $\hat{\varphi}$ is uniformly polynomially unstable if and only if there exist constants $k\in \mathbb{N}$, with $k\geq 2$, and $\theta>1$ such that
\begin{equation}\label{eq.hai.is}
\|\varPhi(k\,n,n,y)w\|\geq \theta \|w\|, \text{ for all $(n,y,w)\in \mathbb{N}^*\times Y\times W$.}
\end{equation}
\end{lem}

\begin{proof}		
{\it Necessity.} Assume that $\hat{\varphi}$ is uniformly polynomially unstable, that is, there  exist $\nu>0$ and $N\geq1$ such that
\begin{equation}\label{upis}
N \|\varPhi(m,n,y)w\| \geq   \left(\frac{m+1}{\;n+1}\right)^{\nu}\| w \|, \; (m,n,y,w)\in \Delta \times Y\times W.		
\end{equation}
	
Letting $k=\left[2N^{\frac{1}{\nu}}\right]$, we get $k>2N^{\frac{1}{\nu}}-1$, thus $\theta:=\frac{1}{N}\left(\frac{k+1}{2}\right)^{\nu}>1$. Then, by \eqref{upis}, we obtain
$$\|\varPhi(k\,n,n,y)w\|\geq \frac{1}{N}\left(\frac{k\,n+1}{n+1}\right)^{\nu} \|w\|\geq \frac{1}{N}\left(\frac{k+1}{2}\right)^{\nu} \|w\|=\theta \| w\|,$$	
for all $(n,y,w)\in \mathbb{N}^*\times Y\times W$, which proves the desired inequality.
	
{\it Sufficiency.} Assume that there exist $k\in\mathbb N$, with $k\geq 2$, and $\theta >1$  such that \eqref{eq.hai.is} holds. Inductively, one can easily check that for any $p\in\mathbb{N}^*$,
\begin{equation}\label{eq.hai.is.p}
\|\varPhi(k^pn,n,y)w\|\geq \theta^p \|w\|, \text{ for all $(n,y,w)\in \mathbb{N}^*\times Y\times W$.}
\end{equation}

Choose $\nu>0$ such that $\theta=k^\nu$ (i.e., $\nu=\frac{\ln{\theta}}{\ln{k}})$, and set $N=N_1M2^{\omega+\nu}$, where $N_1=Mk^{\omega+\nu}$ and $\omega>0$, $M\geq 1$ are the constants given by Definition \ref{defuhis} for $h(n)=n+1$, $n\in\mathbb{N}$.

Let $(m,n,y,w)\in\Delta \times Y\times W$.
If $m>n\geq1$, we consider $p=\left[\dfrac{\ln \frac{m}{n}}{\ln k}\right]$.
Thus, $$n\leq k^p n\leq m<k^{p+1}n.$$

In particular, if $p=0$, then the above inequality becomes $n<m <k\,n$, and the uniform polynomial decay of $\hat{\varphi}$ yields
\begin{align*}
N_1 \| \varPhi(m,n,y)w\|& \geq \frac{N_1}{M} \left( \frac{m+1}{\;n+1} \right)^{-\omega} \| w\|\\
&=\frac{N_1}{M} \left( \frac{m+1}{\;n+1} \right)^{-\omega-\nu} \left( \frac{m+1}{\;n+1} \right)^{\nu} \| w\|\\
&\geq\frac{N_1}{M \, k^{\omega+\nu}} \left( \frac{m+1}{\;n+1} \right)^{\nu} \| w\|\\
&=\left( \frac{m+1}{\;n+1} \right)^{\nu} \| w\|,
\end{align*}
where the last inequality follows from $\frac{m+1}{\;n+1}<k$, as  $m+1\leq kn <kn+k$.

On the other hand, if $p\in\mathbb{N}^*$, then from uniform polynomial decay, inequality \eqref{eq.hai.is.p}, and the choice of $\nu$, we obtain
\begin{align*}
N_1\|\varPhi(m,n,y) w\|
&=N_1\|\varPhi(m,k^p n,\varphi(k^p n,n,y))\varPhi(k^pn,n,y)w\|\\
&\geq\frac{N_1}{M}\left(\frac{k^p n+1}{m+1}\right)^{\omega}\|\varPhi(k^p n,n,y)w\|\\
&\geq\frac{N_1}{M} \left(\frac{k^p n+1}{k^{p+1}n+1}\right)^{\omega}\theta^p \|w\|\\
&\geq \frac{N_1}{M k^\omega}\,k^{\nu p}\|w\|\\
&\geq \frac{N_1}{M k^{\omega+\nu}}\left( \frac{m+1}{\;n+1}\right)^\nu \| w \|\\
&=\left( \frac{m+1}{\;n+1}\right)^\nu \| w \|,
\end{align*}
since
$m+1\leq k^{p+1}n< k^{p+1}(n+1)$, which implies $k^p > \frac{1}{k}\, \frac{m+1}{\;n+1}$. As $\varPhi(n,n,y)=\mathrm{Id}$ and $N_1\geq 1$, this yields
\begin{equation}\label{eq.clar}
N_1\|\varPhi(m,n,y) w\|\geq \left( \frac{m+1}{\;n+1}\right)^\nu \| w \|,\quad  m\geq n\geq1.
\end{equation}

If $m>n=0$, then by \eqref{eq.clar} and \eqref{hbd}, we have
\begin{align*}
N\|\varPhi(m,0,y)w\|&= N_1M2^{\omega+\nu} \| \varPhi(m,1,\varphi(1,0,y))\varPhi (1,0,y)w\|\\
&\geq M2^{\omega+\nu} \left(\frac{m+1}{2}\right)^{\nu}
		\|\varPhi(1,0,y)w\|\\
		&\geq (m+1)^{\nu}\|w\|.
\end{align*}
	
Since $N\geq N_1$, it follows that $\hat{\varphi}$ is uniformly polynomially unstable with exponent $\nu$ and bound $N$.
\end{proof}

By combining the above result with the arguments used in the proof of Proposition \ref{crit log ups}, one readily obtains the following characterization of uniform polynomial instability.

\begin{proposition}\label{crit log upis}\normalfont
	Assume that $\hat{\varphi}$ has uniform polynomial decay. Then   $\hat{\varphi}$ is uniformly polynomially unstable if and only if there exists  a constant $K>0$ such that
	\begin{equation*}\label{eq.log.upis}
	K\|\varPhi(m,n,y)w\|\geq \ln\frac{m+1}{\;n+1}\, \|w\|, \text{ for all $(m,n,y,w)\in \Delta\times Y\times W$.}
	\end{equation*}
\end{proposition}

The following result is the counterpart of Proposition \ref{legatura-uhs-ups-discret} for uniform instability.

\begin{proposition}\label{legatura-uhis-upis-discret}
	Let $h:\mathbb{N}\to [1,\infty)$ be a discrete growth rate satisfying
	\eqref{romb} for some constant $c\geq h(1)$, and assume that $\hat{\varphi}$ has uniform $h$-decay. Then $\hat{\varphi}_{h}$ has uniform polynomial decay.
	
	Moreover,
	$\hat{\varphi}$ is uniformly $h$-unstable if and only if $\hat{\varphi}_{h}$ is uniformly polynomially unstable.		
\end{proposition}
\begin{proof}
Since the proof is similar to that of Proposition \ref{legatura-uhs-ups-discret}, we present only the proof that uniform polynomial instability of $\hat{\varphi}_{h}$ implies uniform $h$-instability of $\hat{\varphi}$. We focus on the case where, for $m,n\in\mathbb{N}$ with $m>n$, there exist $i,j\in\mathbb{N}$ with $i>j$ such that
	$$[\tilde{h}^{-1}(j+1)]\leq n< [\tilde{h}^{-1}(j+2)]\leq [\tilde{h}^{-1}(i+1)]\leq m < [\tilde{h}^{-1}(i+2)].$$

Indeed,  for $(y,w)\in Y\times W$, by successively applying identity $(ii)$ in Definition \ref{def.cocycle}, the uniform  $h$-decay of $\hat{\varphi}$, the uniform polynomial instability of $\hat{\varphi}_{h}$, and relation \eqref{steluta}, and using the inequalities $\frac{i+1}{i+2}\geq \frac{1}{2}$ and $i+1\geq \frac{i+2}{2}$, we obtain
	\begin{align*}\label{eq1.s}
		\| \varPhi(m,n,y)w\|
		&=\|\varPhi(m,[\tilde{h}^{-1}(i+1)],\varphi([\tilde{h}^{-1}(i+1)],n,y)) \\
		&\qquad \varPhi([\tilde{h}^{-1}(i+1)],[\tilde{h}^{-1}(j+2)],\varphi([\tilde{h}^{-1}(j+2)],n,y))\\
		&\qquad \varPhi([\tilde{h}^{-1}(j+2)],n,y)w\|\\
		&\geq \frac{1}{M} \left(\frac{h([\tilde{h}^{-1}(i+1)])}{h(m)}\right)^{\omega}\|\varPhi_h(i,j+1,\varphi([\tilde{h}^{-1}(j+2)],n,y))\\
		&\qquad  \varPhi([\tilde{h}^{-1}(j+2)],n,y)w\|\\
		&\geq \frac{1}{M}\frac{1}{c^\omega} \left(\frac{i+1}{i+2}\right)^{\omega} \frac{1}{N} \left( \frac{i+1}{j+2} \right)^{\nu}
		 \| \varPhi([\tilde{h}^{-1}(j+2)],n,y)w\|\\
		&\geq  \frac{1}{MN (2c)^\omega} \left( \frac{i+1}{j+2} \right)^{\nu}  \frac{1}{M} \left(\frac{h(n)}{h([\tilde{h}^{-1}(j+2)])}\right)^{\omega}\|w\|\\
		&\geq \frac{1}{M^2N (2c)^\omega}   \left( \frac{i+1}{j+2} \right)^{\nu} \frac{1}{c^\omega}\left(\frac{j+1}{j+2}\right)^\omega  \|w\|\\
		&\geq \frac{1}{M^2N (2c)^{2\omega}} \frac{1}{2^\nu} \left(\frac{i+2}{j+2}\right)^{\nu}\|w\|\\
		&\geq \frac{1}{M^2N(2c)^{2\omega+\nu}\,2^\nu} \left(\frac{h(m)}{h(n)}\right)^{\nu}\|w\|,\\
	\end{align*}
where $$\frac{i+2}{j+2}=\frac{h(m)}{h(n)} \, \frac{ h(n) }{j+1}\,\frac{j+1}{j+2}\,\frac{i+2}{h(m)}\geq \frac{1}{2c}\frac{h(m)}{h(n)},$$
since $ \frac{ h(n) }{j+1}\geq \frac{1}{c}$, $\frac{j+1}{j+2}\geq \frac{1}{2}$, and $\frac{i+2}{h(m)}\geq 1$.

In the case where, for $m,n\in\mathbb{N}$ with $m>n$, there exists $i\in\mathbb{N}$ such that
	$$[\tilde{h}^{-1}(i+1)]\leq n<m\leq [\tilde{h}^{-1}(i+2)],$$
we get 	
$$	\| \varPhi(m,n,y)w\| \geq \frac{1}{M\,(2c)^{\omega+\nu}}\left( \frac{h(m)}{h(n)} \right)^{\nu}\|w\|.$$
Therefore $\hat{\varphi}$ is uniformly $h$-unstable.
\end{proof}

To establish the main result of this section, we first need a preliminary result extending Proposition \ref{crit log upis} to uniform $h$-instability. This result provides a  simple characterization of uniform $h$-instability.

\begin{proposition} \label{critlog-uhis-discret}
Let $h:\mathbb{N}\to [1,\infty)$ be a discrete growth rate satisfying
	\eqref{romb} for some constant $c\geq h(1)$.
	If $\hat{\varphi}$ has uniform $h$-decay, then it is uniformly $h$-unstable if and only if there exists a constant $K>0$ such that
	\begin{equation}\label{eq.ln.uhis}
		K	\|\varPhi(m,n,y)w\|\geq \ln \frac{h(m)}{h(n)}\,\|w\|, \text{ for all } (m,n,y,w)\in \Delta \times Y\times W.
	\end{equation}
\end{proposition}

\begin{proof}
	\emph{Necessity.} Using similar arguments as in Proposition \ref{crit log ups}, one can easily prove that \eqref{eq.ln.uhis} holds for $K=\frac{N}{\nu \,e}$, where the constants  $\nu>0$ and $N\geq 1$ are given by Definition \ref{defuhis}.

	\emph{Sufficiency.} Assume that there exists a constant $K>0$ such that \eqref{eq.ln.uhis} holds. Let $(m,n,y,w)\in \Delta \times Y\times W$, and consider $\widetilde{K}=\max\{ (\ln c+1)K, Me^\omega \ln c+K\}$,  where $M\geq 1$ and $\omega>0$ are given by Definition \ref{defuhis}.
	
If $\ln \dfrac{h ([\tilde h^{-1}(m+1)])}{h ([\tilde h^{-1}(\;n+1)] )}\geq 1$, then  \eqref{eq.ln.uhis} and \eqref{eq.nou} imply
	\begin{align*}
\widetilde{K} \|\varPhi_{h}(m,n,y)w\|& \geq (\ln c+1)K	\|\varPhi( [\tilde h^{-1}(m+1)],  [\tilde h^{-1}(n+1)],y)w\|\\
&\geq (\ln c+1) \ln \dfrac{h ([\tilde h^{-1}(m+1)])}{h ([\tilde h^{-1}(\;n+1)] )} \,\|w\|\\
&\geq \ln \frac{m+1}{\;n+1}\, \| w\|.
	\end{align*}
	
On the other hand, if $\ln \dfrac{h ([\tilde h^{-1}(m+1)])}{h ([\tilde h^{-1}(\;n+1)] )}<1$,
	then by \eqref{hbd}, \eqref{eq.ln.uhis}, \eqref{eq.lnmn}, and using the inequality $e > \frac{h ([\tilde h^{-1}(m+1)])}{h ([\tilde h^{-1}(\;n+1)] )}$, we obtain
	\begin{align*}
\widetilde{K} \|\varPhi_{h}(m,n,y)w\|& \geq (Me^\omega \ln c+K)\, \|\varPhi( [\tilde h^{-1}(m+1)],  [\tilde h^{-1}(n+1)],y)w\|\\
&\geq e^\omega \ln c \left( \frac{h ([\tilde h^{-1}(m+1)])}{h ([\tilde h^{-1}(\;n+1)] )} \right)^{-\omega} \|w\|+\ln  \frac{h ([\tilde h^{-1}(m+1)])}{h ([\tilde h^{-1}(\;n+1)] )} \,\| w\|\\
&\geq \left( \ln{c}+ \ln  \frac{h ([\tilde h^{-1}(m+1)])}{h ([\tilde h^{-1}(\;n+1)] )}\right) \| w\|\\
&\geq \ln \frac{m+1}{\;n+1}\, \| w\|.
	\end{align*}

Hence,
$$\widetilde{K} \|\varPhi_{h}(m,n,y)w\|\geq \ln \frac{m+1}{\;n+1}\, \| w\|, \; (m,n,y,w)\in\Delta\times Y\times W,$$
and finally,
Propositions \ref{crit log upis} and \ref{legatura-uhis-upis-discret} complete the proof.
\end{proof}

We can now state the main result of this section, which provides a new criterion of Datko-type for uniform $h$-instability of discrete skew-product evolution semiflows, and in particular for discrete evolution families.

\begin{theorem}\label{t1Zabczykhinst}
Let $h:\mathbb{N}\to [1,\infty)$ be a discrete growth rate such that there exists a constant $c\geq h(1)$ satisfying $$\frac{h(n+1)}{h(n)}\leq c, \text{ for every } n\in \mathbb{N},$$
and assume that $\hat{\varphi}=(\varphi, \varPhi)$ is a discrete skew-product evolution semiflow on $Y\times W$, where the state space $W$ is nontrivial. If  $\hat{\varphi}$  has uniform $h$-decay, then it is uniformly $h$-unstable if and only if
there exists a constant $D>0$ such that
\begin{equation}\label{eq.1inst}
\sum_{k=n}^{\infty}\frac{h(k+1)-h(k)}{h(k)}\,\frac{1}{\|\varPhi(k,n,y)w\|} \leq \frac{D}{\| w\|}, \,\forall\, (n,y,w)\in \mathbb{N} \times Y\times W,\,  w\neq 0.
\end{equation}	
\end{theorem}
\begin{proof}
\emph{Necessity.} Assuming that $\hat{\varphi}$ is uniformly $h$-unstable with constants $\nu>0$ and $N\geq 1$ and using \eqref{eq.halpha}, one can easily show that
\eqref{eq.1inst} holds for $D=\frac{N \, c^{\nu+1}}{\nu}$.

\emph{Sufficiency}. By arguments similar to those used in \emph{Step 1} of Theorem \ref{t1Zabczykhs}, one can deduce that
there exists a constant $K\geq 1$ satisfying  $$K\|\varPhi(m,n,y)w\|\geq  \|w\|,  \text{ for all } (m,n,y,w)\in \Delta \times Y \times W.$$
Then, for $(m,n,y,w)\in\Delta\times Y\times W$ with $m\geq n+1$ and $w\neq 0$, we obtain		
	\begin{align*}
		\frac{1}{\|\varPhi(m,n,y)w\|}\ln \frac{h(m)}{h(n)}&= \sum_{k=n}^{m-1} \frac{1}{\|\varPhi(m,n,y)w\|}\ln \frac{h(k+1)}{h(k)}\\
		&= \sum_{k=n}^{m-1} \frac{1}{\|\varPhi(m,k,\varphi(k,n,y))\varPhi(k,n,y)w\|}\ln \frac{h(k+1)}{h(k)}\\
		&\leq K \sum_{k=n}^{m-1}\frac{1}{\|\varPhi(k,n,y)w\|}\left( \frac{h(k+1)}{h(k)}-1\right)\\
		&\leq K \sum_{k=n}^{\infty}  \frac{h(k+1)-h(k)}{h(k)}\,\frac{1}{\|\varPhi(k,n,y)w\|} \\
		&\leq \frac{K D }{\|w\|},
\end{align*}
thus
$$KD \|\varPhi(m,n,y)w\|\geq \ln \frac{h(m)}{h(n)}\, \| w\|, \text{ for all } (m,n,y,w)\in\Delta\times Y\times W.$$
Proposition \ref{critlog-uhis-discret} now implies that $\hat{\varphi}$ is uniformly $h$-unstable.
\end{proof}

Note that the uniform $h$-decay of $\hat{\varphi}$ implies the injectivity of $\Phi(m,n,y)$, ensuring that the fractions in \eqref{eq.1inst} are well-defined.

The following example shows that the uniform $h$-decay assumption in Theorem \ref{t1Zabczykhinst} cannot be
omitted, even when all operators $\Phi(m,n,y)$ are invertible.

\begin{exemple}\rm
Let $s=(s_k)_{k\in\mathbb{N}}$ be the real sequence considered in Example \ref{ex.g}, and let $h_s$ be its corresponding discrete growth rate defined by \eqref{eq.hgr.s}. The invertible operators
$$ A_n=\begin{cases} \frac{1}{2^j} \, \mathrm{Id},& n=2j\\ 2^{j+1} \, \mathrm{Id},& n=2j+1\end{cases}, \quad j\in\mathbb{N},$$
generate a discrete evolution cocycle $\varPhi:\Delta\times Y\to \mathcal{B}(W)$ over any discrete evolution semiflow $\varphi:\Delta\times Y\to Y$, given by
$$\varPhi(m,n,y)=\begin{cases} A_{m-1} \cdots A_n,& m>n\\
\mathrm{Id},& m=n\end{cases},\quad (m,n,y)\in\Delta\times Y,$$
where $Y$ is an arbitrary base space and $W$ is a nontrivial state space.
Since the operators $\varPhi(m,n,y)$ are invertible, being
a finite product of invertible operators, we get that $\varPhi(m,n,y)w\neq 0$ if $w\neq 0$.

We claim that $\hat{\varphi}=\left(\varphi,\varPhi\right)$ is not uniformly $h_s$-unstable, even though
\begin{equation}\label{eq.1instex}
\sum_{k=n}^{\infty}\frac{h_s(k+1)-h_s(k)}{h_s(k)}\,\frac{1}{\|\varPhi(k,n,y)w\|} \leq \frac{4}{\| w\|},
\end{equation}
for all $(n,y,w)\in \mathbb{N} \times Y\times W$ with  $w\neq 0$.
\end{exemple}

\begin{proof}[\unskip\nopunct]
Indeed, if we assume that $\hat{\varphi}$ is uniformly $h_s$-unstable, then there exist constants $\nu>0$ and $N\geq 1$ such that
$$N \|\varPhi(m,n,y)w\| \geq   \left(\frac{h_s(m)}{h_s(n)}\right)^{\nu}\| w \|,$$
for all $(m,n,y,w)\in \Delta \times Y\times W.$
Taking $m=2k+1$, $n=0$,  and $w\neq 0$, and using the fact that  $\varPhi(2k+1,0,y)=\mathrm{Id}$, we obtain
$$ N\geq \left(h_s(2k+1)\right)^\nu,\quad k\in\mathbb{N},$$
which is false, since $h_s(2k+1)\to\infty$ as $k\to\infty$. Therefore,  $\hat{\varphi}$ is not uniformly $h_s$-unstable.

Let $(n,y,w)\in \mathbb{N} \times Y\times W$ with $w\neq 0$. To prove \eqref{eq.1instex}, we consider two cases.

\medskip

{\it Case 1}. $n=2j$, $j\in\mathbb{N}$.

\medskip

Notice first that
$$\varPhi(2i,2j,y)=2^{i-j}\,\mathrm{Id} \text{ and } \varPhi(2i+1,2j,y)=\frac{1}{2^j}\,\mathrm{Id}, \; \forall\, i\geq j.$$
This yields
\begin{align*}
\sum_{k=n}^{\infty}\frac{h_s(k+1)-h_s(k)}{h_s(k)}\,\frac{1}{\|\varPhi(k,n,y)w\|}&=\sum_{k=2j}^{\infty}\frac{s_k}{\|\varPhi(k,2j,y)w\|}\\
&=\sum_{i=j}^{\infty}\left( \frac{1}{2^{i-j}}+\frac{1}{2^i}\, 2^j\right)\frac{1}{\| w \|}\\
&=2\left( \sum_{k=0}^{\infty}\frac{1}{2^k}\right)\frac{1}{\| w \|}\\
&=\frac{4}{\| w \|}.
\end{align*}

\medskip

{\it Case 2}. $n=2j+1$, $j\in\mathbb{N}$.

\medskip

In this case, let us observe that
$$\varPhi(k,2j+1,y)=\varPhi(k,2j+2,y)A_{2j+1}=2^{j+1}\varPhi(k,2j+2,y),\quad k\geq 2j+2.$$
Then, using the definition of $s_k$ and the previous case, we obtain
\begin{align*}
\sum_{k=n}^{\infty}\frac{h_s(k+1)-h_s(k)}{h_s(k)}\,\frac{1}{\|\varPhi(k,n,y)w\|}
&=\sum_{k=2j+1}^{\infty}\frac{s_k}{\|\varPhi(k,2j+1,y)w\|}\\
&=\frac{1}{2^j}\frac{1}{\| w \|}+\frac{1}{2^{j+1}}\sum_{k=2j+2}^{\infty}\frac{s_k}{\|\varPhi(k,2j+2,y)w\|}\\
&=\left( \frac{1}{2^j} +\frac{4}{2^{j+1}}\right)\frac{1}{\| w \|}\\
&=\frac{3}{2^j}\frac{1}{\| w \|}.
\end{align*}
This proves \eqref{eq.1instex}.
\end{proof}

Setting $h(n)=e^n$ and $h(n)=n+1$  in Theorem \ref{t1Zabczykhinst},  we obtain the following discrete Datko-type characterizations of uniform exponential instability and uniform polynomial instability, respectively, for discrete skew-product evolution semiflows with nontrivial state spaces.

\begin{corollary}\label{cexpis}
	A discrete skew-product evolution semiflow $\hat{\varphi}=(\varphi, \varPhi)$ with uniform exponential decay is uniformly exponentially unstable if and only if there exists a constant $D>0$ such that
	\begin{equation*}\label{eq.exp.is}
		\sum_{k=n}^{\infty} \frac{1}{\|\varPhi(k,n,y)w\|} \leq \frac{D}{ \|w\|}, \text{ for all } (n,y,w)\in  \mathbb{N}\times Y\times W \text{ with }w\neq 0.
	\end{equation*}
\end{corollary}

\begin{corollary}\label{cpolis}
	A discrete skew-product evolution semiflow $\hat{\varphi}=(\varphi, \varPhi)$ with uniform polynomial decay is uniformly polynomially unstable if and only if there exists a constant $D>0$ such that
	\begin{equation*}\label{eq.pol.is}
		\sum_{k=n}^{\infty} \frac{1}{(k+1)\|\varPhi(k,n,y)w\|} \leq \frac{D}{ \|w\|}, \text{ for all } (n,y,w)\in  \mathbb{N}\times Y\times W \text{ with }w\neq 0.
	\end{equation*}
\end{corollary}

\section{Some examples}

In what follows, we provide a concrete example of a skew-product evolution semiflow satisfying the assumptions of Theorem \ref{t1Zabczykhs}.

\begin{exemple}\rm\label{ex1Ts}
Let $Y=[1,2]^{\mathbb{N}}$ be the set of all sequences $(y_n)_{n\in\mathbb{N}}$ taking values in the interval $[1,2]$.
The mapping $\varphi:\Delta\times Y\to Y$, defined by
$$\varphi(m,n,y)=L^{m-n}y,$$
where $L:Y\to Y$ is the left shift operator on $Y$, i.e., $L(y_0,y_1,\ldots):=(y_1,y_2,\ldots),$
is a discrete  evolution semiflow on $Y$.

Assume that $W$ is an arbitrary Banach space.
For $(m,n,y)\in \Delta\times Y$, we define $$\varPhi(m,n,y)=
\begin{cases}
\prod\limits_{k=n}^{m-1} \left(1+s_k y_{k-n}\right)^{-1}\mathrm{Id},& m>n,\\
\mathrm{Id},& m=n,
\end{cases}
$$
where $s=(s_k)_{k\in\mathbb{N}}$ is a sequence of positive real numbers satisfying the assumptions in Example \ref{ex.gr}. Then, $\hat{\varphi}=\left(\varphi,\varPhi\right)$ is a discrete skew-product evolution semiflow on $Y\times W$. Indeed, for $m>n>n_0$ in $\mathbb{N}$ and $y=(y_0,y_1,\ldots)\in Y$, we have
\begin{align*}
\varPhi(m,n,\varphi(n,n_0,y))&\varPhi(n,n_0,y)\\
&=\prod\limits_{k=n}^{m-1} \left(1+s_k \left(L^{n-n_0}y\right)_{k-n}\right)^{-1}
\prod\limits_{k=n_0}^{n-1} \left(1+s_k y_{k-n_0}\right)^{-1}\mathrm{Id}\\
&=\prod\limits_{k=n}^{m-1} \left(1+s_k y_{k-n_0}\right)^{-1}
\prod\limits_{k=n_0}^{n-1} \left(1+s_k y_{k-n_0}\right)^{-1}\mathrm{Id}\\
&=\prod\limits_{k=n_0}^{m-1} \left(1+s_k y_{k-n_0}\right)^{-1}\mathrm{Id}\\
&=\varPhi(m,n_0,y),
\end{align*}
where we have used the fact that
$L^{n-n_0}y=\left(y_{n-n_0},y_{n-n_0+1},\ldots\right),$
and hence
$$\left(L^{n-n_0}y\right)_{k-n}=y_{n-n_0+k-n}=y_{k-n_0}.$$

We claim that $\hat{\varphi}$ satisfies the assumptions of Theorem \ref{t1Zabczykhs} with respect to the discrete growth rate $h_s$ defined in Example \ref{ex.gr}, thus  it is uniformly $h_s$-stable.
\end{exemple}

\begin{proof}[\unskip\nopunct]
Evidently, $\hat{\varphi}$ is uniformly $h_s$-bounded since
$$\| \varPhi(m,n,y)\|\leq 1\leq \frac{h_s(m)}{h_s(n)},\quad (m,n,y)\in\Delta\times Y.$$

Fix $(n,y,w)\in \mathbb{N}\times Y\times W$ and let
$$a_k:=\| \varPhi(k,n,y)w \|,\quad k\geq n.$$
Then
\begin{equation*}
a_{k+1}=\prod\limits_{j=n}^{k} \left(1+s_j y_{j-n}\right)^{-1}\| w\|=\frac{1}{1+s_k y_{k-n}}\, a_k.
\end{equation*}
This yields
\begin{align*}
a_k-a_{k+1}&= \left(1-\frac{1}{1+s_k y_{k-n}}\right) a_k
= \frac{s_k y_{k-n}}{1+s_k y_{k-n}} a_k\geq \frac{s_k}{1+s_k} a_k,
\end{align*}
since $s_k>0$, $y_{k-n}\geq 1$, and the function $f_k(x)=\frac{s_k x}{1+s_k x}$ is nondecreasing on $[1,\infty)$. Furthermore, using the inequality $s_k\leq S$, we obtain
$$a_k-a_{k+1}\geq \frac{s_k}{1+S} a_k,$$
and hence
$$s_k a_k\leq (1+S)\left(a_k-a_{k+1}\right).$$
Finally, for $m>n$, we have
\begin{align*}
\sum_{k=n}^{m-1}\frac{h_s(k+1)-h_s(k)}{h_s(k)}\,\|\varPhi(k,n,y)w\|&=\sum_{k=n}^{m-1} s_k a_k\\
&\leq (1+S)\sum_{k=n}^{m-1}\left(a_k-a_{k+1}\right)\\
&=(1+S)(a_n-a_m)\\
&\leq (1+S)\|w\|.
\end{align*}
Letting $m\to\infty$, we conclude that \eqref{eq.1} holds with $D=1+S$. Therefore, by Theorem \ref{t1Zabczykhs},  $\hat{\varphi}$ is uniformly $h_s$-stable.
\end{proof}

The following example illustrates the applicability of Theorem \ref{t1Zabczykhinst}.

\begin{exemple}\rm
Let $W$ be an arbitrary nontrivial Banach space. Using the same base space $Y$, discrete evolution semiflow $\varphi$, and sequence $s=(s_k)_{k\in\mathbb{N}}$ as in Example \ref{ex1Ts}, we define
$$\varPsi(m,n,y)=
\begin{cases}
\prod\limits_{k=n}^{m-1} \left(1+s_k y_{k-n}\right)\mathrm{Id},& m>n,\\
\mathrm{Id},& m=n.
\end{cases}
$$
One can prove that the pair  $(\varphi,\varPsi)$ is a discrete skew-product evolution semiflow on $Y\times W$ satisfying
$$\| \varPsi (m,n,y)w \|\geq \| w\| \geq \left(\frac{h_s(m)}{h_s(n)}\right)^{-1} \|w\|,\quad (m,n,y,w)\in \Delta\times Y\times W,$$
and hence it has uniform $h_s$-decay. On the other hand, for $w\neq 0$, we have
$$\frac{1}{\| \varPsi (k,n,y)w \|}=\frac{\| \varPhi (k,n,y)w \|}{\|w\|^2},$$
where $\varPhi$ is the discrete evolution cocycle in Example \ref{ex1Ts}. Consequently,
$$\sum_{k=n}^{\infty}\frac{h_s(k+1)-h_s(k)}{h_s(k)}\,\frac{1}{\|\varPsi(k,n,y)w\|} \leq \frac{1+S}{\| w\|},$$
for all $(n,y,w)\in\mathbb{N}\times Y\times W$ with $w\neq 0$.
Therefore, by Theorem \ref{t1Zabczykhinst},  $(\varphi,\varPsi)$ is uniformly $h_s$-unstable.
\end{exemple}

%
%
%

\end{document}